\documentclass{amsart}

\usepackage{amsmath,amsthm,amssymb,latexsym,amscd}
\usepackage{graphicx,psfrag}

\newcommand{\C}{\mathbb{C}}
\newcommand{\Z}{\mathbb{Z}}
\newcommand{\comment}[1]{}

\newtheorem{thm}{Theorem}
\newtheorem{lem}[thm]{Lemma}
\newtheorem{coro}[thm]{Corollary}
\newtheorem{propo}[thm]{Proposition}
\newtheorem{defi}[thm]{Definition}

\newtheorem{rem}[thm]{Remark}
\newtheorem{nota}[thm]{Notation}
\newtheorem{ques}[thm]{Question}

\begin{document}

\psfrag{a}{$\!\alpha$}
\psfrag{a-r}{$\!\alpha-r$}
\psfrag{-b}{$\;\beta$}
\psfrag{-b-s}{$\;\beta-s$}
\psfrag{c}{$\!\gamma$}
\psfrag{c-t}{$\!\gamma-t$}

\psfrag{i}{$i$}
\psfrag{i-1}{$\!\!\!i-1$}
\psfrag{j}{$j$}
\psfrag{j-1}{$j-1$}
\psfrag{-j}{$-j$}
\psfrag{j-b}{$j-b$}
\psfrag{k}{$k$}
\psfrag{m}{$m$}
\psfrag{m+n}{$m+n$}
\psfrag{n}{$n$}
\psfrag{-n}{$\!-n$}
\psfrag{n-1}{$n-1$}
\psfrag{n-2}{$n-2$}
\psfrag{n-3}{$n-3$}
\psfrag{n+1}{$n+1$}
\psfrag{n+2}{$n+2$}
\psfrag{n+3}{$n+3$}
\psfrag{p}{$p$}
\psfrag{r}{$r$}
\psfrag{s}{$s$}
\psfrag{t}{$t$}
\psfrag{...}{$\cdots$}
\psfrag{. . .}{$\cdots$}
\psfrag{-2}{$-2$}
\psfrag{-1}{$-1$}
\psfrag{1}{$1$}
\psfrag{2}{$2$}
\psfrag{3}{$3$}

\title[skein module of surgery on a $(2,2b)$ torus link]{On the Kauffman bracket skein module of surgery on a $(2,2b)$ torus link}

\author{John M. Harris}
\address{University of Southern Mississippi\\
Long Beach, Mississippi}
\email{john.m.harris@usm.edu}
\comment{
}

\begin{abstract} 
We show that the Kauffman bracket skein modules of certain manifolds obtained from integral surgery on a $(2,2b)$ torus link are finitely generated, and list the generators for select examples. 
\end{abstract}

\maketitle

\section{Introduction}
\label{sec:Introduction}

In \cite{Ka88}, Kauffman presents an elegant construction of the Jones polynomial, an invariant of oriented links in $S^3$, by constructing a new invariant, the Kauffman bracket polynomial.  The Kauffman bracket is an invariant of unoriented framed links in $S^3$, defined by the following skein relations:

\begin{enumerate}
	\item $\left<\begin{minipage}{0.4in}\includegraphics[width=0.4in]{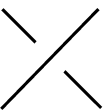}\end{minipage}\right> = A \left<\;\begin{minipage}{0.4in}\includegraphics[width=0.4in]{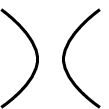}\end{minipage}\;\right> + A^{-1} \left<\;\begin{minipage}{0.4in}\includegraphics[width=0.4in]{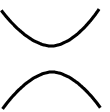}\end{minipage}\;\right>$
	\item $<L \cup \hbox{unknot}> = (-A^{-2} - A^2) <L>$
\end{enumerate}
 
For the invariant to be well-defined, one also must normalize it by choosing a value for the empty link.  $<\hbox{empty link}> = 1$, for instance.

Alternatively, we can use the skein relations to construct a module of equivalence classes of links in $S^3$, or, for that matter, in any oriented 3-manifold.  See Przytycki \cite{Pr91} and Turaev \cite{Tu88}.

\begin{defi}
Let $N$ be an oriented 3-manifold, and let $R$ be a commutative ring with identity, with a specified unit $A$.  The Kauffman bracket skein module of $N$, denoted $S(N;R,A)$, or simply $S(N)$, is the free $R$-module generated by the framed isotopy classes of unoriented links in $N$, including the empty link, quotiented by the skein relations which define the Kauffman bracket.
\end{defi}

Since every crossing and unknot can be eliminated from a link in $S^3$ by the skein relations, $S(S^3)$ is generated by the empty link. Kauffman's argument that his bracket polynomial is well-defined shows that $S(S^3)$ is free on the empty link.

For $R = \Z[A^{\pm 1}]$, Hoste and Przytycki have computed the skein modules of all of the closed, oriented manifolds of genus 1:  $S(L(p,q))$ , which is free on $\left\lfloor\frac{p}{2}\right\rfloor + 1$ generators \cite{HP93}, and $S(S^1 \times S^2) \cong \Z[A^{\pm 1}]\oplus
(\bigoplus^\infty_{i=1}\Z[A^{\pm 1}]/(1-A^{2i+4}))$ \cite{HP95}.  Over $\Z[A^{\pm 1}]$, localized by inverting all of the cyclotomic polynomials, Gilmer and the author have computed the skein module of the quaternionic manifold \cite{GH}.

Additionally, Bullock has determined whether or not the skein module obtained from integral surgery on a trefoil is finitely generated in \cite{Bu97}.  In this paper, we pursue a similar result, for integral surgery on a $(2,2b)$ torus link.

\begin{figure}

$$\begin{minipage}{0.4in}\includegraphics{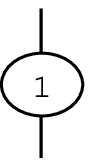}\end{minipage} 
= \;\;\begin{minipage}{0.45in}\includegraphics{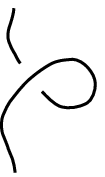}\end{minipage}
= \;\;\begin{minipage}{0.45in}\includegraphics[angle=180]{postwist1r}\end{minipage}
\;\;\;\;\;\;\;\;
\begin{minipage}{0.55in}\includegraphics{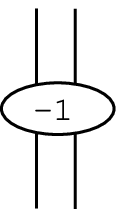}\end{minipage} 
= \;\;\begin{minipage}{0.75in}\includegraphics{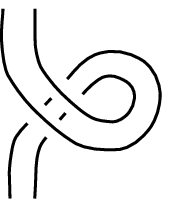}\end{minipage}
= \;\;\begin{minipage}{0.75in}\includegraphics[angle=180]{negtwist2r}\end{minipage}$$

\caption{Examples of twist notation}
\label{fig:twistexamples}
\end{figure}

\begin{nota}
For any integer $n$, let $\;\begin{minipage}{0.7in}\includegraphics{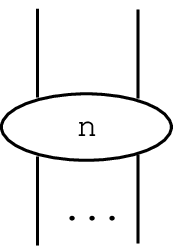}\end{minipage}$ denote $n$ full twists in the depicted strands.  For example, see Figure \ref{fig:twistexamples}.  
\end{nota}

\begin{defi} 
We define $M(\alpha,\beta,\gamma)$ to be the manifold obtained by surgery on the torus link 
\;\;$\begin{minipage}{2.4in}\includegraphics{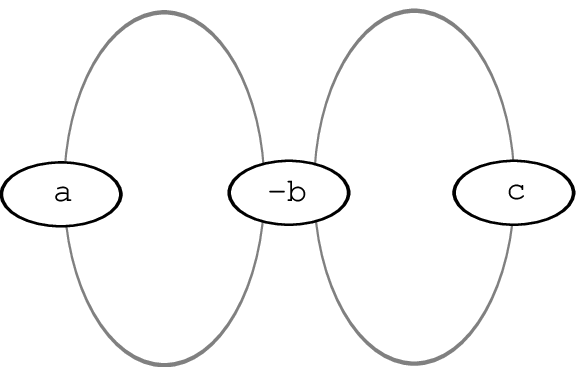}\end{minipage}$ , with the blackboard framing.
\end{defi}

\begin{thm}
For all integers $\alpha$, $\beta$, and $\gamma$ such that $a = |\alpha| > 1,b = |\beta| > 1,c = |\gamma| > 1$, $\frac{1}{a} < \frac{1}{b} + \frac{1}{c}$, $\frac{1}{b} < \frac{1}{a} + \frac{1}{c}$, and $\frac{1}{c} < \frac{1}{a} + \frac{1}{b}$, $S(M(\alpha,\beta,\gamma))$ is finitely generated.
\label{thm:main thm}
\end{thm}

For specific values of $\alpha$, $\beta$, and $\gamma$, we can use brute-force computation to refine our result, explicitly listing generating sets for $S(M(\alpha,\beta,\gamma))$.

\begin{nota}
We refer to the collection of loops 
$$\begin{minipage}{2.4in}\includegraphics{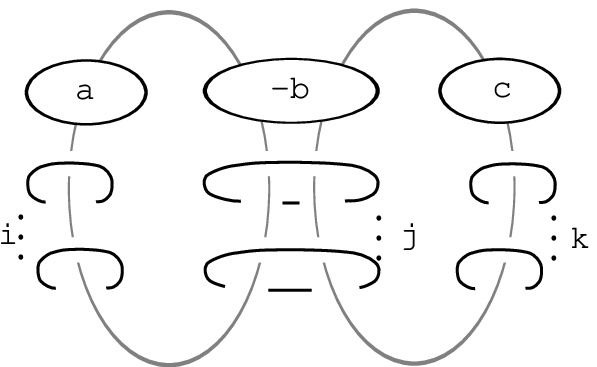}\end{minipage}$$
in $M(\alpha,\beta,\gamma)$ using the algebraic notation $x^i y^j z^k$.
\end{nota}

In particular, we obtain the following for $S(M(2,-2,2))$, $S(M(3,-2,3))$, and $S(M(3,-2,5))$ (the skein modules of the 3-fold, 4-fold, and 5-fold branched cyclic coverings of $S^3$ over the trefoil, respectively, as listed by Rolfsen \cite{Ro76}): 

$$\begin{array}{ccccc}
\underline{\alpha}&\underline{\beta}&\underline{\gamma}&\mbox{ \underline{fundamental group} }&\mbox{ \underline{generators} }\\
\\
2&-2&2&\mbox{ quaternion group }& 1, z, z^2, y, x\\
\\
3&-2&3&\mbox{ binary tetrahedral group }& 1, z, z^2, z^3, y, x, x^2\\
\\
3&-2&5&\;\;\;\;\;\;\;\;\mbox{ binary icosahedral group }\;\;\;\;\;\;\;\;& 1, z, z^2, z^3, z^4, z^5, y, x, x^2
\end{array}$$

Note that the generating set for the skein module of the quaternionic manifold essentially coincides with what was shown in \cite{GH} over the ring $R'$ obtained from $\Z[A^{\pm 1}]$ by inverting the multiplicative set generated by the elements of the set $\{A^n - 1|n \in \Z^+\}$.  Since any dependence relation over $\Z[A^{\pm 1}]$ would hold over $R'$ and since $S(M(2,-2,2);R',A)$ is a free module of rank 5, we obtain the following:

\begin{coro}
$S(M(2,-2,2);\Z[A^{\pm 1}],A)$ is a free module of rank 5.
\end{coro}

\section{Twists and Loops}
\label{sec:Twists and Loops}

\begin{figure}

$$\begin{minipage}{0.6in}\includegraphics{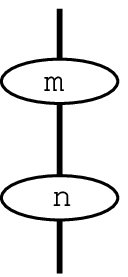}\end{minipage} 
= \;\;\;\begin{minipage}{1in}\includegraphics{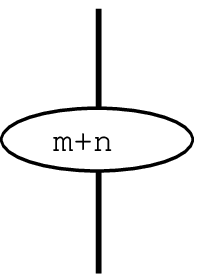}\end{minipage}
\;\;\;\;\;\;\;
\begin{minipage}{1.2in}\includegraphics{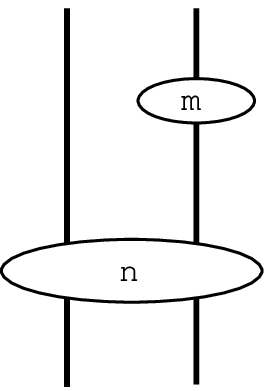}\end{minipage} 
= \;\;\;\begin{minipage}{1in}\includegraphics{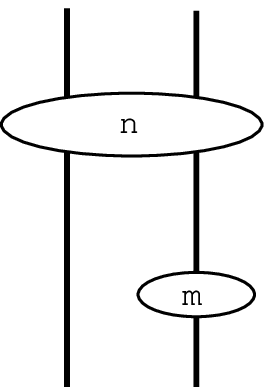}\end{minipage}$$

$$\begin{minipage}{1.25in}\includegraphics{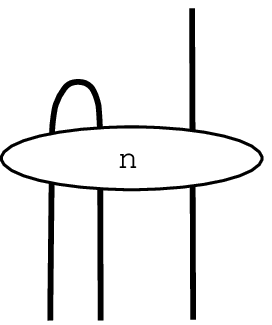}\end{minipage} 
= \;\;\;\begin{minipage}{1.25in}\includegraphics{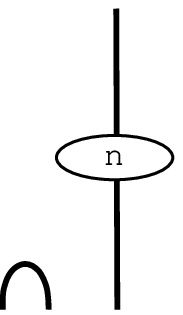}\end{minipage}$$

\caption{Useful properties of twists}
\label{fig:twistproperties}
\end{figure}

Twists have many useful properties, a few of which are listed in Figure \ref{fig:twistproperties}.  Note that, to obtain clearer diagrams, we represent a fixed but arbitrary number of parallel strands with a thick line.  

\begin{figure}

$$\begin{minipage}{1in}\includegraphics{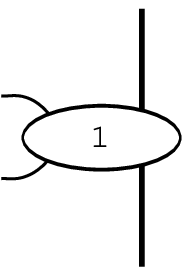}\end{minipage} 
= \;\;\begin{minipage}{1.9in}\includegraphics{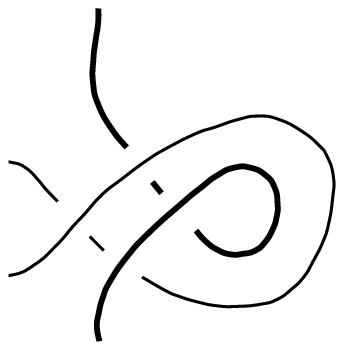}\end{minipage} 
= \;\;\begin{minipage}{1.9in}\includegraphics{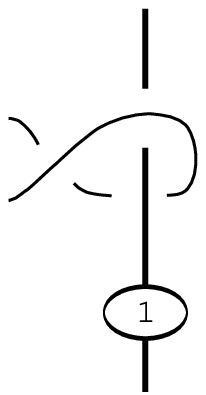}\end{minipage}$$

$$\begin{minipage}{1in}\includegraphics{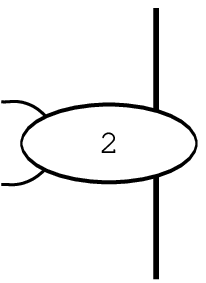}\end{minipage} 
= \;\;\;\begin{minipage}{1.25in}\includegraphics{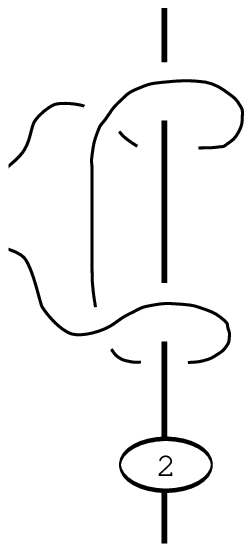}\end{minipage}$$ 

$$= A \begin{minipage}{1.25in}\includegraphics{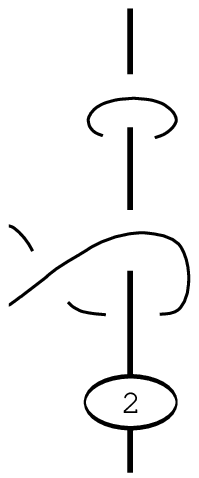}\end{minipage} 
+ A^{-1} \begin{minipage}{1.25in}\includegraphics{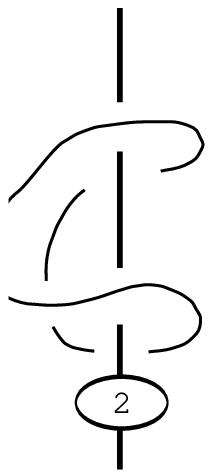}\end{minipage}$$

$$= A \begin{minipage}{1.25in}\includegraphics{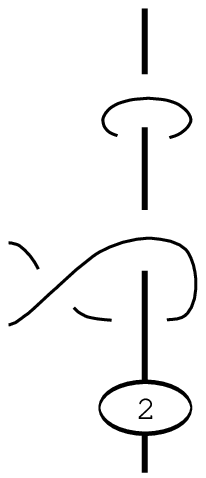}\end{minipage} 
- A^2 \;\;\;\begin{minipage}{1.25in}\includegraphics{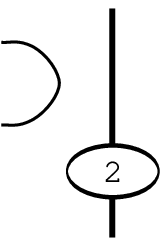}\end{minipage}$$

\caption{Examples of rewriting twists}
\label{fig:opentwistsasloops}
\end{figure}

We are most interested in using skein relations and isotopy to rewrite one strand, twisted with others, as a linear combination involving loops encircling the others, as in Figure \ref{fig:opentwistsasloops}.

In fact, repeating by repeating the steps performed in Figure \ref{fig:opentwistsasloops}, we obtain the following lemma:

\begin{lem}  For each integer $n > 0$, 

$$\begin{minipage}{0.9in}\includegraphics{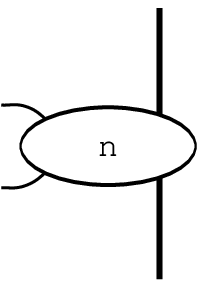}\end{minipage} 
= \sum_{j<n} f^+_j(A) \begin{minipage}{1.25in}\includegraphics{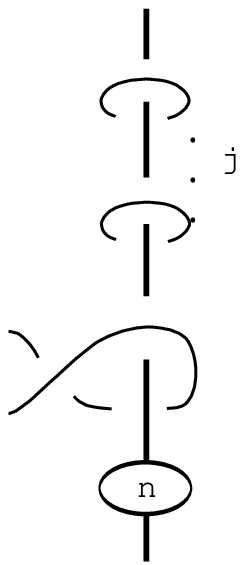}\end{minipage}
+ \sum_{j<n-1} g^+_j(A) \;\;\;\;\begin{minipage}{1in}\includegraphics{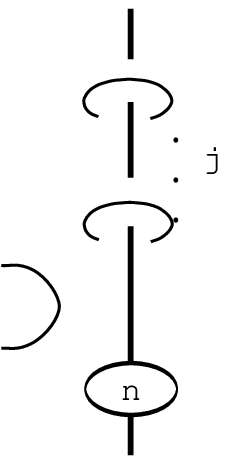}\end{minipage},$$

where $f^+_{n-1}(A) = A^{n-1}$, and 

$$\begin{minipage}{0.9in}\includegraphics{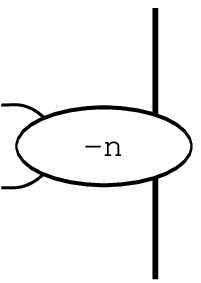}\end{minipage} 
= \sum_{j<n} f^-_j(A) \begin{minipage}{1.25in}\includegraphics{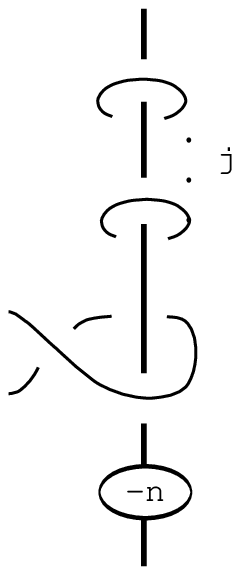}\end{minipage}
+ \sum_{j<n-1} g^-_j(A) \;\;\;\;\begin{minipage}{1in}\includegraphics{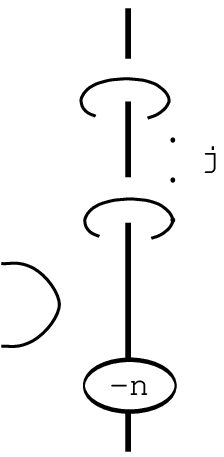}\end{minipage},$$

where $f^-_{n-1}(A) = A^{1-n}$.

\label{lem:opentwist1}
\end{lem}

\begin{proof} For $n = 1$ and $n=2$, the result is obtained in Figure \ref{fig:opentwistsasloops}.  

Let $n > 2$, and suppose that the result holds for all $k < n$.  Then

$$\begin{minipage}{1in}\includegraphics{nopentwistl}\end{minipage} 
= \;\;\;\;\begin{minipage}{1.25in}\includegraphics{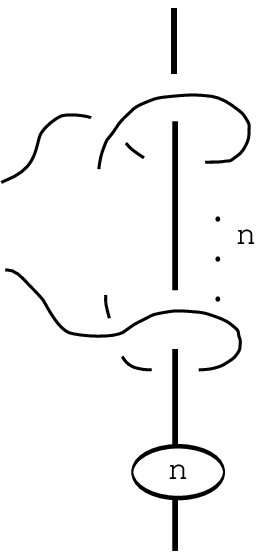}\end{minipage}
= \;\;\;\;\begin{minipage}{1.25in}\includegraphics{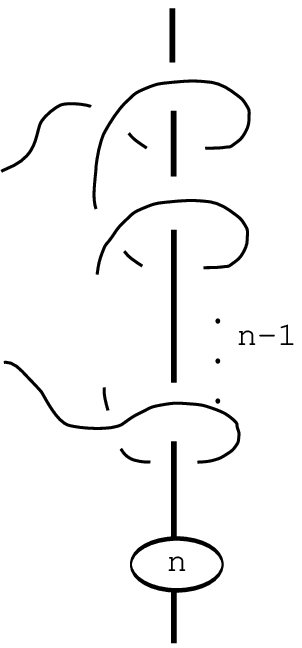}\end{minipage}$$

$$= A \;\;\;\;\begin{minipage}{1.5in}\includegraphics{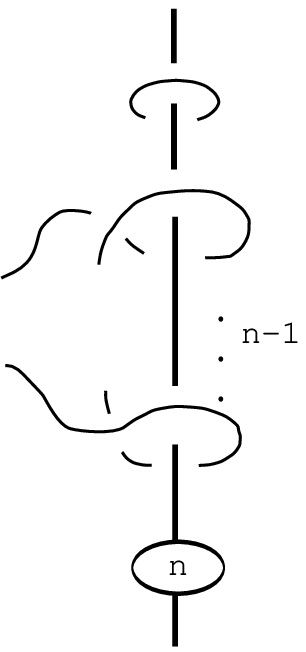}\end{minipage} 
+ A^{-1} \;\;\;\;\begin{minipage}{1.5in}\includegraphics{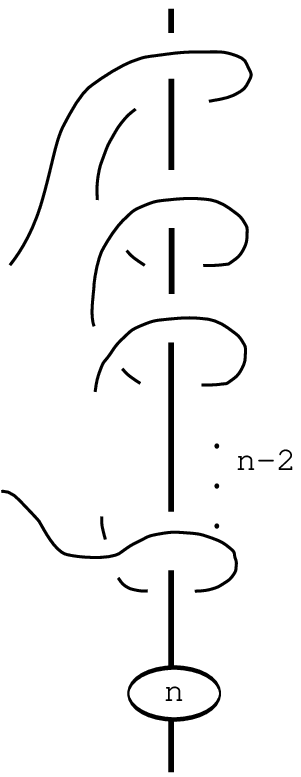}\end{minipage}$$

$$= A \;\;\;\;\begin{minipage}{1in}\includegraphics{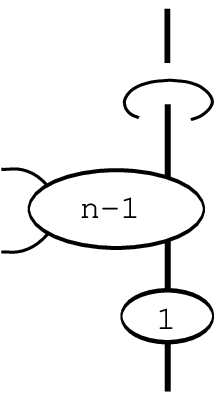}\end{minipage} 
- A^2 \;\;\;\;\begin{minipage}{1in}\includegraphics{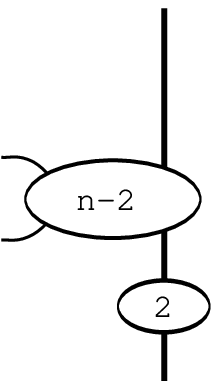}\end{minipage}$$

Hence, the first equation follows by induction on $n$.  The second equation can be obtained by reversing all of the crossings in the first.  \end{proof}

By rotating the diagrams in the previous lemma by 180 degrees, we also obtain

\begin{lem}  For each integer $n > 0$, 

$$\begin{minipage}{1in}\includegraphics{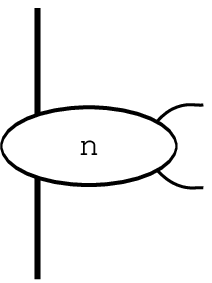}\end{minipage} 
= \sum_{j<n} f^+_j(A) \begin{minipage}{1.25in}\includegraphics{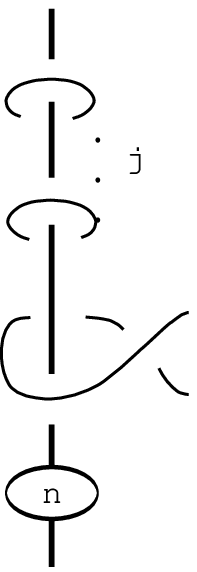}\end{minipage}
+ \sum_{j<n-1} g^+_j(A) \;\;\;\;\begin{minipage}{1in}\includegraphics{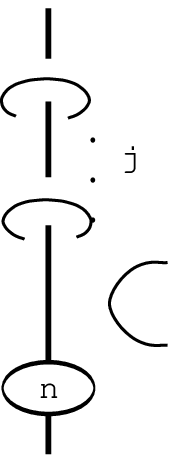}\end{minipage},$$

where $f^+_{n-1}(A) = A^{n-1}$, and 

$$\begin{minipage}{1in}\includegraphics{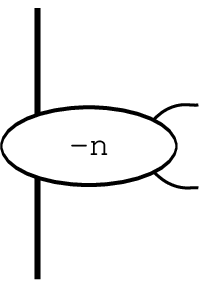}\end{minipage} 
= \sum_{j<n} f^-_j(A) \begin{minipage}{1.25in}\includegraphics{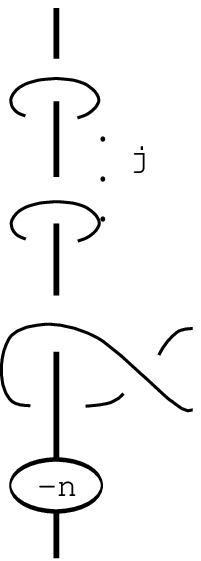}\end{minipage}
+ \sum_{j<n-1} g^-_j(A) \;\;\;\;\begin{minipage}{1in}\includegraphics{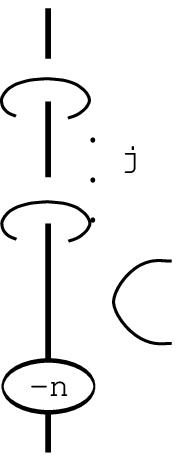}\end{minipage},$$

where $f^-_{n-1}(A) = A^{1-n}$.

\label{lem:opentwist2}
\end{lem}

In particular, if a component of a link is only twisted about one set of other strands, we obtain, as an immediate corollary of Lemma \ref{lem:opentwist1}, 

\begin{lem} For each integer $n > 0$,

$$\begin{minipage}{1in}\includegraphics{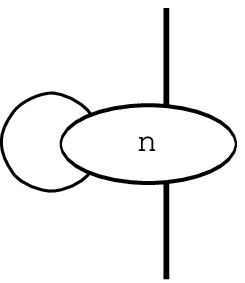}\end{minipage}
= \sum_{i \leq n} h^+_i(A) \;\;\;\;\begin{minipage}{.75in}\includegraphics{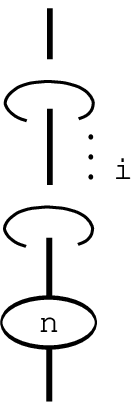}\end{minipage},$$

where $h^+_n(A) = -A^{n+2}$, and

\psfrag{n}{$\!\!\!-n$}
$$\begin{minipage}{1in}\includegraphics{nclosedtwist}\end{minipage}
= \sum_{i \leq n} h^-_i(A) \;\;\;\;\begin{minipage}{.75in}\includegraphics{iloops}\end{minipage},$$
\psfrag{n}{$n$}

where $h^-_n(A) = -A^{-n-2}$.

\label{lem:closedtwist1}
\end{lem}

Similarly, as a corollary of Lemma \ref{lem:opentwist2},

\begin{lem} For each integer $n > 0$,

$$\begin{minipage}{1in}\includegraphics{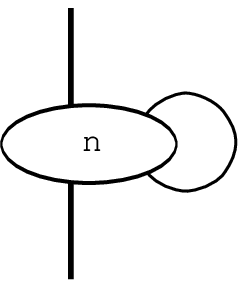}\end{minipage}
= \sum_{i \leq n} h^+_n(A) \;\;\;\;\begin{minipage}{.75in}\includegraphics{iloops}\end{minipage},$$

where $h^+_n(A) = -A^{n+2}$, and

\psfrag{n}{$\!\!\!-n$}
$$\begin{minipage}{1in}\includegraphics{revnclosedtwist}\end{minipage}
= \sum_{i \leq n} h^-_n(A) \;\;\;\;\begin{minipage}{.75in}\includegraphics{iloops}\end{minipage},$$
\psfrag{n}{$n$}

where $h^-_n(A) = -A^{-n-2}$.

\label{lem:closedtwist2}
\end{lem}

Suppose that a component of a link is twisted with two sets of strands.  While more complicated than in the cases previously considered, it is still possible to rewrite the component as a linear combination of loops around the other strands:

\psfrag{r}{$\!m$}
\psfrag{s}{$\!n$}
\begin{lem}  For all integers $m, n > 0$,

$$\begin{minipage}{1.3in}\includegraphics{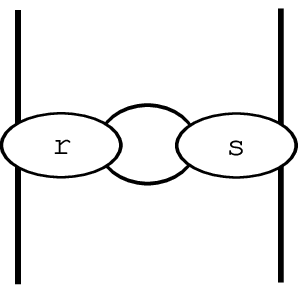}\end{minipage} 
= \sum_{i \leq m, j \leq n} f^{++}_{i,j}(A) \;\;\;\;
\begin{minipage}{2in}\includegraphics{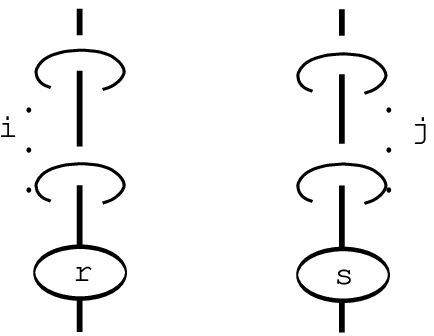}\end{minipage}$$ 

$$+ \sum_{i < m, j < n} g^{++}_{i,j}(A) \;\;\;\;
\begin{minipage}{2in}\;\;\;\;\includegraphics{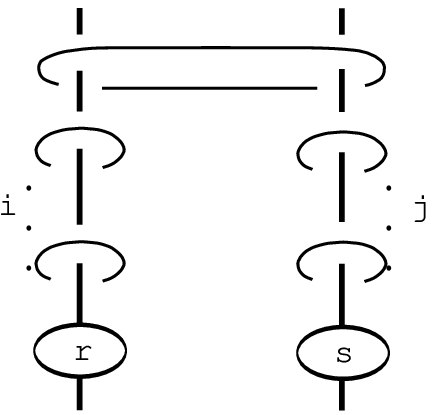}\end{minipage},$$

where $f^{++}_{m,n}(A) = -A^{m+n+2}$.

\label{lem:dtwistpp}
\end{lem}

\begin{proof} Applying Lemma \ref{lem:opentwist1}, and then applying Lemma \ref{lem:opentwist2} to each diagram of the resulting linear combination, 

$$\begin{minipage}{1.3in}\includegraphics{dtwist}\end{minipage} 
= \sum_{i < m, j < n} f^+_i(A) f^+_j(A) \;\;\;\;
\begin{minipage}{2in}\includegraphics{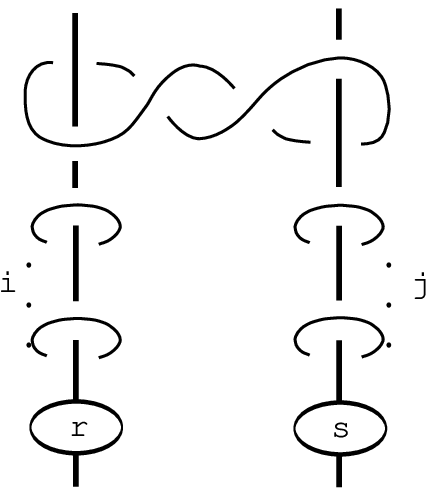}\end{minipage}$$ 

$$+ \sum_{i \leq m, j < n-1} (-A^3) f^+_i(A) g^+_j(A) \;\;\;\;
\begin{minipage}{2in}\;\;\;\;\includegraphics{dtwistr1}\end{minipage}$$

$$+ \sum_{i < m-1, j \leq n} (-A^3) g^+_i(A) f^+_j(A) \;\;\;\;
\begin{minipage}{2in}\;\;\;\;\includegraphics{dtwistr1}\end{minipage}$$

$$+ \sum_{i < m-1, j < n-1} (-A^2 - A^{-2}) g^+_i(A) g^+_j(A) \;\;\;\;
\begin{minipage}{2in}\;\;\;\;\includegraphics{dtwistr1}\end{minipage}.$$

Since $\;\;\;\;\begin{minipage}{1.3in}\includegraphics{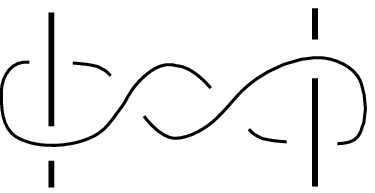}\end{minipage}$

$$ = -A^4 \;\;\;\;\begin{minipage}{1.67in}\includegraphics{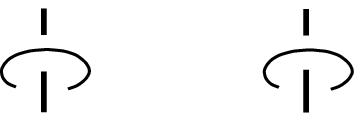}\end{minipage}
- A^2 \;\;\;\;\begin{minipage}{1.67in}\includegraphics{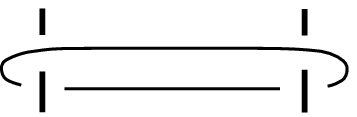}\end{minipage},$$

the result follows.

\end{proof}

\psfrag{r}{$\!m$}
\psfrag{s}{$\!\!\!-n$}
\begin{lem}  For all integers $m, n > 0$,

$$\begin{minipage}{1.3in}\includegraphics{dtwist}\end{minipage} 
= \sum_{\scriptsize \begin{array}{c}i \leq m, j < n-1\\i < m-1, j \leq n\end{array}} f^{+-}_{i,j}(A) \;\;\;\;
\begin{minipage}{2in}\includegraphics{dtwistr1}\end{minipage}$$ 

$$+ \sum_{i < m, j < n} g^{+-}_{i,j}(A) \;\;\;\;
\begin{minipage}{2in}\;\;\;\;\includegraphics{dtwistr2}\end{minipage}.$$

where $g^{+-}_{m,n}(A) = A^{m-n}$.

\label{lem:dtwistpn}
\end{lem}

\begin{proof}  Applying Lemma \ref{lem:opentwist1}, and then applying Lemma \ref{lem:opentwist2} to each diagram of the resulting linear combination, 

$$\begin{minipage}{1.3in}\includegraphics{dtwist}\end{minipage} 
= \sum_{i < m, j < n} f^+_i(A) f^-_j(A) \;\;\;\;
\begin{minipage}{2in}\includegraphics{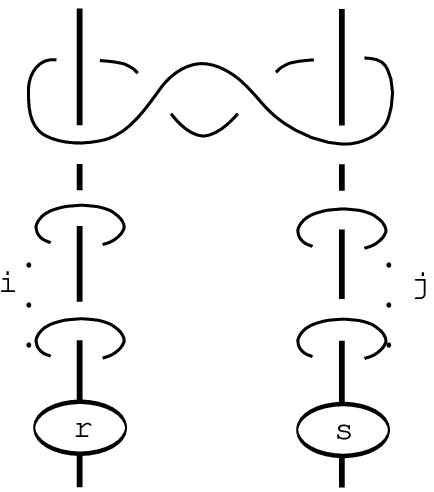}\end{minipage}$$ 

$$+ \sum_{i \leq m, j < n-1} (-A^3) f^+_i(A) g^-_j(A) \;\;\;\;
\begin{minipage}{2in}\;\;\;\;\includegraphics{dtwistr1}\end{minipage}$$

$$+ \sum_{i < m-1, j \leq n} (-A^{-3}) g^+_i f^-_j(A) \;\;\;\;
\begin{minipage}{2in}\;\;\;\;\includegraphics{dtwistr1}\end{minipage}$$

$$+ \sum_{i < m-1, j < n-1} (-A^2 - A^{-2}) g^+_i(A) g^-_j(A) \;\;\;\;
\begin{minipage}{2in}\;\;\;\;\includegraphics{dtwistr1}\end{minipage}.$$

\end{proof}

By an argument similar to that for  Lemma \ref{lem:dtwistpn}, we obtain 

\psfrag{r}{$\!\!\!\!-m$}
\psfrag{s}{$\!n$}
\begin{lem}  For all integers $m, n > 0$,

$$\begin{minipage}{1.3in}\includegraphics{dtwist}\end{minipage} 
= \sum_{\scriptsize \begin{array}{c}i \leq m, j < n-1\\i < m-1, j \leq n\end{array}} f^{-+}_{i,j}(A) \;\;\;\;
\begin{minipage}{2in}\includegraphics{dtwistr1}\end{minipage}$$ 

$$+ \sum_{i < m, j < n} g^{-+}_{i,j}(A) \;\;\;\;
\begin{minipage}{2in}\;\;\;\;\includegraphics{dtwistr2}\end{minipage},$$

where $g^{-+}_{m,n}(A) = A^{n-m}$.

\label{lem:dtwistnp}
\end{lem}

\psfrag{r}{$\!\!\!\!-m$}
\psfrag{s}{$\!\!\!-n$}
By an argument similar to that for Lemma \ref{lem:dtwistpp}, we obtain 
\begin{lem}  For all integers $m, n > 0$,

$$\begin{minipage}{1.3in}\includegraphics{dtwist}\end{minipage} 
= \sum_{i \leq m, j \leq n} f^{--}_{i,j}(A) \;\;\;\;
\begin{minipage}{2in}\includegraphics{dtwistr1}\end{minipage}$$ 

$$+ \sum_{i < m, j < n} g^{--}_{i,j}(A) \;\;\;\;
\begin{minipage}{2in}\;\;\;\;\includegraphics{dtwistr2}\end{minipage},$$

where $f^{--}_{i,j}(A) = -A^{-m-n-2}$.

\label{lem:dtwistnn}
\end{lem}

\section{Finitely Generating the Skein Module}
\label{sec:Finitely Generating the Skein Module}

\psfrag{r}{$r$}
\psfrag{s}{$s$}

Since all links in the exterior of the surgery description of $M(\alpha,\beta,\gamma)$ can be isotoped into a genus 2 handlebody and since the skein relations allow us to remove all crossings in a diagram, $S(M(\alpha,\beta,\gamma))$ is generated by $\{x^i y^j z^k\}$.

\begin{defi} For $a = |\alpha|,b=|\beta|,c=|\gamma| > 0$, we define a strict linear ordering on the generating set $\{x^i y^j z^k\}$ of $M(\alpha,\beta,\gamma)$ as follows:
$x^i y^j z^k < x^m y^n z^p$ if
\begin{itemize}
\item $\frac{i}{a} + \frac{j}{b} + \frac{k}{c} < \frac{m}{a} + \frac{n}{b} + \frac{p}{c}$,
\item $\frac{i}{a} + \frac{j}{b} + \frac{k}{c} = \frac{m}{a} + \frac{n}{b} + \frac{p}{c}$ and $i(k+1) < m(p+1)$,
\item $\frac{i}{a} + \frac{j}{b} + \frac{k}{c} = \frac{m}{a} + \frac{n}{b} + \frac{p}{c}$, $i(k+1) = m(p+1)$, and
$\max\left(\frac{j}{b},\frac{k}{c}\right) < \max\left(\frac{n}{b},\frac{p}{c}\right)$,
\item $\frac{i}{a} + \frac{j}{b} + \frac{k}{c} = \frac{m}{a} + \frac{n}{b} + \frac{p}{c}$, $i(k+1) = m(p+1)$, 
$\max\left(\frac{j}{b},\frac{k}{c}\right) = \max\left(\frac{n}{b},\frac{p}{c}\right)$, and $j < n$, or
\item $\frac{i}{a} + \frac{j}{b} + \frac{k}{c} = \frac{m}{a} + \frac{n}{b} + \frac{p}{c}$, $i(k+1) = m(p+1)$, 
$\max\left(\frac{j}{b},\frac{k}{c}\right) = \max\left(\frac{n}{b},\frac{p}{c}\right)$, $j = n$, and $k < p$.
\end{itemize}
\end{defi}

Suppose that $a,b,c > 1$, $\frac{1}{a} < \frac{1}{b} + \frac{1}{c}$, $\frac{1}{b} < \frac{1}{a} + \frac{1}{c}$, and $\frac{1}{c} < \frac{1}{a} + \frac{1}{b}$.

By sliding over an attached 2-handle, we obtain 

\begin{defi} the Type I relation: 

$$\begin{minipage}{2.3in}\includegraphics{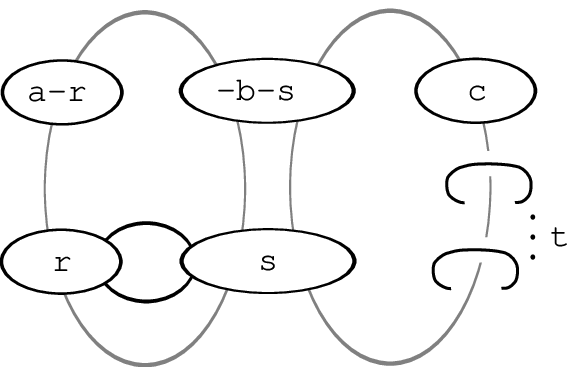}\end{minipage}  
= \begin{minipage}{2.3in}\includegraphics{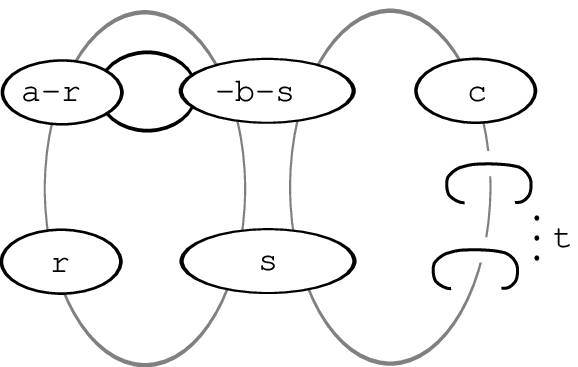}\end{minipage}$$
\end{defi}
First, note that by Lemmas \ref{lem:dtwistpp}, \ref{lem:dtwistpn}, \ref{lem:dtwistnp}, and \ref{lem:dtwistnn}, each side of the relation can be written as a linear combination of loops of the form $x^i y^j z^k$, since for all nonnegative integers $u, v,$ and $w$, 

\psfrag{i}{$\!u$}
\psfrag{j}{$v$}
\psfrag{k}{$w$}
$$\begin{minipage}{2.1in}\includegraphics[scale=.9]{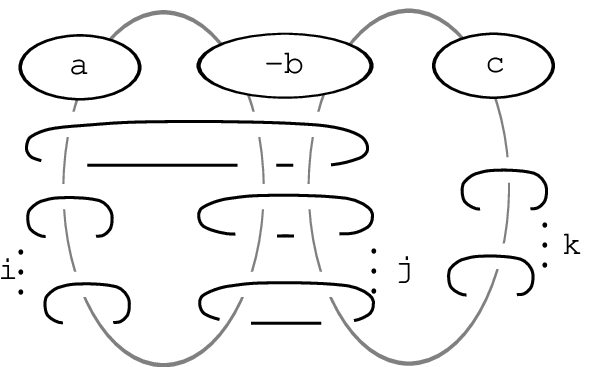}\end{minipage}
\psfrag{k}{$w+1$}
= \;\begin{minipage}{2.1in}\includegraphics[scale=.9]{loops}\end{minipage}\;\;\;\;\;\;\;\;.$$
\psfrag{i}{$i$}
\psfrag{j}{$j$}
\psfrag{k}{$k$}

Note that when $r \geq 0$ and $s \geq 0$, the greatest term appearing on the left side of the Type I relation, rewritten as a linear combination of loops, is $x^r y^s z^t$:

When $r,s > 0$, by Lemma \ref{lem:dtwistpp}, $x^r y^s z^t$ and $x^{r-1} y^{s-1} z^{t+1}$ appear as the greatest terms of their respective types.

Since $\frac{1}{c} < \frac{1}{a} + \frac{1}{b}$, 
$$\frac{r}{a} + \frac{s}{b} + \frac{t}{c} > \left(\frac{r}{a} + \frac{s}{b} + \frac{t}{c}\right) + \left(-\frac{1}{a} - \frac{1}{b} + \frac{1}{c}\right) = \frac{r-1}{a} + \frac{s-1}{b} + \frac{t+1}{c}.$$

When either $r = 0$ or $s = 0$, the claim follows by Lemma \ref{lem:closedtwist1} or Lemma \ref{lem:closedtwist2}.  When both are $0$, the claim follows trivially.

Also note that as long as $r > 0$ or $s > 0$, the leading coefficient is $-A^{r+s+2}$.

Similarly, when $r \leq 0$ and $s \leq 0$, the greatest term appearing on the left side of the Type I relation is $x^{-r} y^{-s} z^t$, and as long as both are not $0$, its coefficient is $-A^{r+s-2}$.

When $r > 0$ and $s < 0$, the greatest term appearing on the left side of the Type I relation is $x^{r-1} y^{-s-1} z^{t+1}$: 

By Lemma \ref{lem:dtwistpn}, $x^{r-1} y^{-s-1} z^{t+1}$, $x^{r-2} y^{-s} z^t$, and $x^r y^{-s-2} z^t$ appear as the greatest terms of their respective types.  Since $\frac{1}{b} < \frac{1}{a} + \frac{1}{c}$, 
$$\frac{r-1}{a} + \frac{-s-1}{b} + \frac{t+1}{c} > \left(\frac{r-1}{a} + \frac{-s-1}{b} + \frac{t+1}{c}\right) + \left(-\frac{1}{a} + \frac{1}{b} - \frac{1}{c}\right) = x^{r-2} y^{-s} z^t.$$  Since $\frac{1}{a} < \frac{1}{b} + \frac{1}{c}$, 
$$\frac{r-1}{a} + \frac{-s-1}{b} + \frac{t+1}{c} > \left(\frac{r-1}{a} + \frac{-s-1}{b} + \frac{t+1}{c}\right) + \left(\frac{1}{a} - \frac{1}{b} - \frac{1}{c}\right) = x^{r} y^{-s-2} z^t.$$  

Also note that in this case, the leading coefficient is $A^{r+s}$.

Similarly, when $r < 0$ and $s > 0$, the greatest term appearing on the left side is $x^{-r-1} y^{s-1} z^{t+1}$, with coefficient $A^{r+s}$.

Likewise, the greatest term on the right side is $x^{|\alpha-r|-1} y^{|\beta-s|-1} z^{t+1}$, when $\alpha-r$ and $\beta-s$ are nonzero with different signs, and the greatest term on the right side is $x^{|\alpha-r|} y^{|\beta-s|} z^t$ otherwise.

By sliding over the other attached 2-handle, we obtain 

\begin{defi} the Type II relation:

$$\begin{minipage}{2.3in}\includegraphics{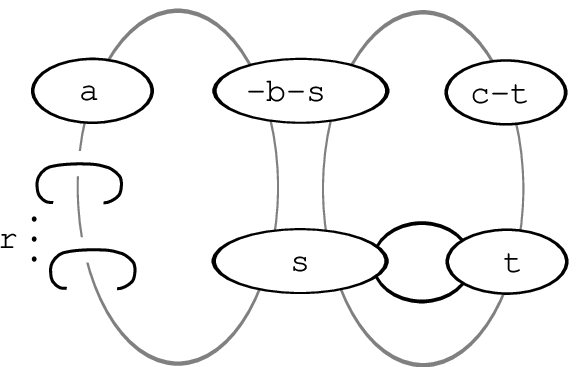}\end{minipage} 
= \begin{minipage}{2.3in}\includegraphics{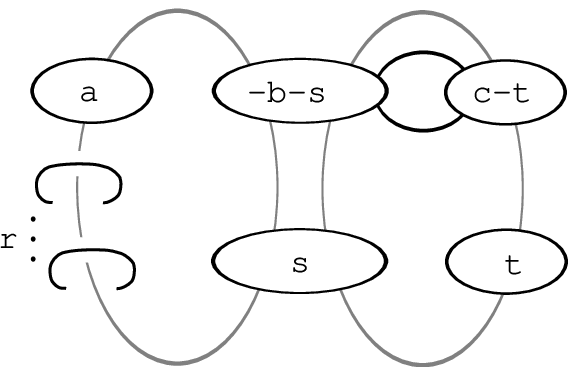}\end{minipage}$$
\end{defi}

As with the Type I relation, each side of the relation can be rewritten as a linear combination of loops of the form $x^i y^j z^k$.

Also, as with the Type I relation, the greatest term appearing on the left side of the Type II relation is $x^{r+1} y^{|s|-1} z^{|t|-1}$ when the signs of $s$ and $t$ differ, with coefficient $A^{s+t}$.  Otherwise, the greatest term appearing on the left side is $x^r y^{|s|} z^{|t|}$, and as long as one of $s$ and $t$ are nonzero, the leading coefficient is $-A^{s+t \pm 2}$.

Finally, as with the Type I relation, the greatest term on the right side of the Type II relation is $x^{r+1} y^{|\beta-s|-1} z^{|\gamma-t|-1}$ when the signs of $\beta-s$ and $\gamma-t$ differ, and the greatest term appearing on the left side is $x^r y^{|\beta-s|} z^{|\gamma-t|}$ otherwise.

\begin{thm} For all integers $a,b,c > 1$ such that $\frac{1}{a} < \frac{1}{b} + \frac{1}{c}$, $\frac{1}{b} < \frac{1}{a} + \frac{1}{c}$, and $\frac{1}{c} < \frac{1}{a} + \frac{1}{b}$, $S(M(a,b,c))$ is finitely generated.
\label{thm:main thm same}
\end{thm}

\begin{proof}  We show that with respect to our previously defined ordering, $x^i y^j z^k$ can be rewritten as linear combinations of lesser terms whenever $i \geq a$, $j \geq b$, or $k \geq c$.  We accomplish this by choosing a Type I or Type II relation in which $x^i y^j z^k$ appears as the greatest term on the left side, as in the previous discussion.  We then show that $x^i y^j z^k$ is greater than the greatest term on the right side of the relation.  Hence, by subtracting all of the terms less than $x^i y^j z^k$ from both sides of the equation and dividing both sides by the (invertible, as previously discussed) coefficient of $x^i y^j z^k$, we successfully rewrite $x^i y^j z^k$.

Case 1:  Suppose $i \geq a$.  Let $r  = i, s = j, t = k$.  Since $r > 0$ and $s \geq 0$, $x^i y^j z^k$ is the greatest term on the left of the Type I relation.  Since $a-r = a-i \leq 0$, the greatest term on the right side is $x^{i-a} y^{j-b} z^k$, if $j \geq b$ or $i = a$, and $x^{i-a-1} y^{b-j-1} z^{k+1}$, if $j < b$ and $i > a$.

Case 1.1:  Suppose $j \geq b$ or $i = a$.  $\frac{i}{a} + \frac{j}{b} + \frac{k}{c} > \frac{i-a}{a} + \frac{j-b}{b} + \frac{k}{c}$, and thus, $x^i y^j z^k > x^{i-a} y^{j-b} z^k$.

Case 1.2:  Suppose $j < b$ and $i > a$.  $\frac{i}{a} + \frac{j}{b} + \frac{k}{c} > \frac{i}{a} - \frac{j}{b} + \frac{k}{c} = \frac{i-a}{a} + \frac{b-j}{b} + \frac{k}{c} > \left(\frac{i-a}{a} + \frac{b-j}{b} + \frac{k}{c}\right) + \left(-\frac{1}{a} - \frac{1}{b} + \frac{1}{c}\right) = \frac{i-a-1}{a} + \frac{b-j-1}{b} + \frac{k+1}{c}$.  Hence, $x^i y^j z^k > x^{i-a-1} y^{b-j-1} z^{k+1}$.

Case 2:  Suppose $i < a$ and $j \geq b$.  Let $r  = i, s = j, t = k$.  Since $r \geq 0$ and $s > 0$, $x^i y^j z^k$ is the greatest term on the left of the Type I relation.  Since $a-r = a-i > 0$ and $b-s = b-j \leq 0$, the greatest term on the right side is $x^{a-i-1} y^{j-b-1} z^{k+1}$, if $j > b$, and $x^{a-i} z^k$, if $j = b$.

Case 2.1:  Suppose $j > b$.  $\frac{i}{a} + \frac{j}{b} + \frac{k}{c} > -\frac{i}{a} + \frac{j}{b} + \frac{k}{c} = \frac{a-i}{a} + \frac{j-b}{b} + \frac{k}{c} > \left(\frac{a-i}{a} + \frac{j-b}{b} + \frac{k}{c}\right) + \left(-\frac{1}{a} - \frac{1}{b} + \frac{1}{c}\right) = \frac{a-i-1}{a} + \frac{j-b-1}{b} + \frac{k+1}{c}$, and thus, $x^i y^j z^k > x^{a-i-1} y^{j-b-1} z^{k+1}$.

Case 2.2:  Suppose $j = b$.  $\frac{i}{a} + \frac{j}{b} + \frac{k}{c} = \frac{i}{a} + 1 + \frac{k}{c} > -\frac{i}{a} + 1 + \frac{k}{c}  = \frac{a-i}{a} + \frac{k}{c}$, and hence, $x^i y^j z^k > x^{a-i} z^k$.

Case 3:  Suppose $i < a$, $j < b$, and $k \geq c$.  Let $r  = i, s = j, t = k$.  Since $s \geq 0$ and $t > 0$, $x^i y^j z^k$ is the greatest term on the left of the Type II relation.  Since $c-t = c-k \leq 0$, the greatest term on the right side is $x^{i+1} y^{b-j-1} z^{k-c-1}$, if $k > c$, and $x^i y^{b-j}$, if $k = c$.

Case 3.1:  Suppose $k > c$.  Then $\frac{i}{a} + \frac{j}{b} + \frac{k}{c} > \frac{i}{a} - \frac{j}{b} + \frac{k}{c} = \frac{i}{a} + \frac{b-j}{b} + \frac{k-c}{c} > \left(\frac{i}{a} + \frac{b-j}{b} + \frac{k-c}{c}\right) + \left(\frac{1}{a} - \frac{1}{b} - \frac{1}{c}\right) = \frac{i+1}{a} + \frac{b-j-1}{b} + \frac{k-c-1}{c}$, and thus, $x^i y^j z^k > x^{i+1} y^{b-j-1} z^{k-c-1}$.

Case 3.2:  Suppose $k = c$.  $\frac{i}{a} + \frac{j}{b} + \frac{k}{c} = \frac{i}{a} + \frac{j}{b} + 1 > \frac{i}{a} - \frac{j}{b} + 1  = \frac{i}{a} + \frac{b-j}{b}$, and so, $x^i y^j z^k > x^i y^{b-j}$.

\end{proof}

\begin{rem}
Note that we can refine the generating set obtained in the above proof, through additional applications of the Type I and Type II relations.  For instance, we can rewrite $x^i y^j z^k$ when 

\begin{itemize}
\item $i < a$, $j < b$, and $\frac{i}{a} + \frac{j}{b} > 1$,
\item $i < a$, $j < b$, $\frac{i}{a} + \frac{j}{b} = 1$, and $i > \frac{a}{2}$,
\item $j < b$, $k < c$ and $\frac{j}{b} + \frac{k}{c} > 1$, or
\item $j < b$, $k < c$, $\frac{j}{b} + \frac{k}{c} = 1$, and $k > \frac{c}{2}$.
\end{itemize}
\end{rem}

\begin{thm} For all integers $a,b,c > 1$ such that $\frac{1}{a} < \frac{1}{b} + \frac{1}{c}$, $\frac{1}{b} < \frac{1}{a} + \frac{1}{c}$, and $\frac{1}{c} < \frac{1}{a} + \frac{1}{b}$, $S(M(a,-b,c))$ is finitely generated.
\label{thm:main thm diff}
\end{thm}

\begin{proof}  We show that with respect to our previously defined ordering, $x^i y^j z^k$ can be rewritten as linear combinations of lesser terms whenever $i \geq a$, $j \geq b$, or $k > c(2-\frac{2}{b})$.  As in the previous proof, we accomplish this by choosing a Type I or Type II relation in which $x^i y^j z^k$ appears as the greatest term on the left side, and then show that $x^i y^j z^k$ is greater than the greatest term on the right side of the relation.  Here, however, the task is a bit more difficult:  the difference in signs prevents us from proceeding in a completely straightforward manner.

Case 1:  Suppose $i \geq a$.  Let $r  = i, s = j, t = k$.  Since $r > 0$ and $s \geq 0$, $x^i y^j z^k$ is the greatest term on the left of the Type I relation.  $a-r = a-i \leq 0$ and $-b-s = -b-j < 0$, and thus $x^{i-a} y^{b+j} z^k$ is the greatest term on the right.  $\frac{i}{a} + \frac{j}{b} + \frac{k}{c} = \frac{i-a}{a} + \frac{b+j}{b} + \frac{k}{c}$, and $i(k+1) > (i-a)(k+1)$, so $x^i y^j z^k > x^{i-a} y^{b+j} z^k$.

Case 2:  Suppose $i < a$ and $j \geq b$.

Case 2.1:  Suppose $k > 0$.  Let $r  = i+1, s = -j-1, t = k-1$.  Since $r > 0$ and $s < 0$, $x^i y^j z^k = x^{(i+1)-1} y^{-(-j-1)-1} z^{(k-1)+1}$ is the greatest term on the left of the Type I relation.  Since $a-r = a-i-1 \geq 0$ and $-b-s = -b+j+1 > 0$, $x^{a-i-1} y^{-b+j+1} z^{k-1}$ is the greatest term on the right.  $\frac{i}{a} + \frac{j}{b} + \frac{k}{c} > \left(\frac{i}{a} + \frac{j}{b} + \frac{k}{c}\right) + \left(-\frac{1}{a} + \frac{1}{b} - \frac{1}{c}\right) = \frac{a-i-1}{a} + \frac{-b+j+1}{b} + \frac{k-1}{c}$, and thus, $x^i y^j z^k > x^{a-i-1} y^{-b+j+1} z^{k-1}$.

Case 2.2:  Suppose $k = 0$.

Case 2.2.1:  Suppose $i > 0$.  Let $r  = i-1, s = -j-1, t = 1$.  Since $s < 0$ and $t > 0$, $x^i y^j$ is the greatest term on the left of the Type II relation.  $-b-s = -b+j+1 > 0$ and $c-t = c-1 > 0$, and thus $x^{i-1} y^{-b+j+1} z^{c-1}$ is the greatest term on the right.  $\frac{i}{a} + \frac{j}{b} + \frac{k}{c} > \left(\frac{i}{a} + \frac{j}{b} + \frac{k}{c}\right) + \left(-\frac{1}{a} + \frac{1}{b} - \frac{1}{c}\right) = \frac{i-1}{a} + \frac{-b+j+1}{b} + \frac{c-1}{c}$, and thus, $x^i y^j > x^{i-1} y^{-b+j+1} z^{c-1}$.

Case 2.2.2:  Suppose $i = 0$.  Let $r  = 0, s = -j, t = 0$.  Since $t = 0$, $y^j$ is the greatest term on the left of the Type II relation.  $-b-s = -b+j \geq 0$ and $c-t = c > 0$, and thus $y^{-b+j} z^c$ is the greatest term on the right.  $\frac{j}{b} = \frac{-b+j}{b} + \frac{c}{c}$ and $0(0+1) = 0(c+1)$.  When $j > b$, $\max\left(\frac{j}{b}, 0\right) > \max\left(\frac{-b+j}{b}, \frac{c}{c}\right)$, and when $j = b$, $\max\left(\frac{j}{b}, 0\right) = 1 = \max\left(\frac{-b+j}{b}, \frac{c}{c}\right)$ and $j = b > 0 = -b+j$.  Hence, $y^j > y^{-b+j} z^c$.

Case 3:  Suppose $i < a, j < b,$ and $k > c(2-\frac{2}{b})$.  (Hence, $k > c$.)

Case 3.1:  Suppose $i > 0$.  Let $r  = i-1, s = -j-1, t = k+1$.  Since $s < 0$ and $t > 0$, $x^i y^j z^k$ is the greatest term on the left of the Type II relation.  Since $-b-s = -b+j+1 \leq 0$ and $c-t = c-k-1 < 0$, $x^{i-1} y^{b-j-1} z^{k-c+1}$ is the greatest term on the right.  $\frac{i}{a} + \frac{j}{b} + \frac{k}{c} > \left(\frac{i}{a} + \frac{j}{b} + \frac{k}{c}\right) + \left(-\frac{1}{a} - \frac{1}{b} + \frac{1}{c}\right) = \frac{i-1}{a} + \frac{b-j-1}{b} + \frac{k-c+1}{c}$, and thus, $x^i y^j > x^{i-1} y^{b-j-1} z^{k-c+1}$.

Case 3.2:  Suppose $i = 0$.

Case 3.2.1:  Suppose $j = b-1$.  Let $r = 1, s = -b, t = k-1$.  Since $r > 0$ and $s < 0$, $y^{b-1} z^k$ is the greatest term on the left of the Type I relation.  $a-r = a-1 > 0$ and $-b-s = 0$, and thus, $x^{a-1} z^{k-1}$ is the greatest term on the right.  Since $\frac{b-1}{b} + \frac{k}{c} > \left(\frac{b-1}{b} + \frac{k}{c}\right) + \left(-\frac{1}{a} + \frac{1}{b} - \frac{1}{c}\right) = \frac{a-1}{a} + \frac{k-1}{c}$,  $y^{b-1} z^k > x^{a-1} z^{k-1}$.

Case 3.2.2:  Suppose $j < b-1$.  Let $r = 0, s = j, t = k$.  Since $s \geq 0$ and $t > 0$, $y^j z^k$ is the greatest term on the left of the Type II relation.  $-b-s = -b-j < 0$ and $c-k < 0$, and thus, $y^{b+j} z^{k-c}$ is the greatest term on the right.  $\frac{j}{b} + \frac{k}{c} = \frac{b+j}{b} + \frac{k-c}{c}$, $0(k+1) = 0(k-c+1)$, and  $\max\left(\frac{j}{b}, \frac{k}{c}\right) = \frac{k}{c} > \max\left(\frac{b+j}{b}, \frac{k-c}{c}\right)$ since $k > c\left(\frac{2b-2}{b}\right) \geq c\left(\frac{b+j}{b}\right)$.  Hence $y^j z^k > y^{b+j} z^{k-c}$.  \end{proof}

\begin{proof}[Proof of Theorem \ref{thm:main thm}]
If $\alpha$, $\beta$, and $\gamma$ are all positive, the result follows by Theorem \ref{thm:main thm same}.  If $\alpha$, $\beta$, and $\gamma$ are all negative, the result follows as well, since $S(M(\alpha,\beta,\gamma))$ is isomorphic to $S(M(-\alpha,-\beta,-\gamma))$.

Suppose that exactly one of $\alpha$, $\beta$, and $\gamma$ is negative.  If $\beta < 0$, the result follows by Theorem \ref{thm:main thm diff}.  If $\alpha < 0$, by sliding the right handle over the left and performing isotopy, we see that $M(\alpha,\beta,\gamma)$ is identical to $M(\gamma,\alpha,\beta)$, and so the result follows.  Similarly, if $\gamma < 0$, by sliding the left handle over the right, $M(\alpha,\beta,\gamma)$ is seen to be identical to $M(\beta,\gamma,\alpha)$, and so again the result follows.

If exactly one of $\alpha$, $\beta$, and $\gamma$ is positive, $S(M(-\alpha,-\beta,-\gamma))$ is finitely generated, and thus $S(M(\alpha,\beta,\gamma))$ is finitely generated as well.
\end{proof}

\section{Examples}
\label{sec:Examples}

While the previous proofs yield a finite set of generators for $S(M(\alpha,\beta,\gamma))$, they do not exploit the full potential of the Type I and Type II relations.  Using the following Python code, we can refine our results for $S(M(a,-b,c))$:
\\

\begin{verbatim}
def greaterthan(a,b,c,i,j,k,m,n,p):
    if i*b*c + j*a*c + k*a*b > m*b*c + n*a*c + p*a*b:
        return True
    elif i*b*c + j*a*c + k*a*b == m*b*c + n*a*c + p*a*b:
        if i*(k+1) > m*(p+1):
            return True
        elif i*(k+1) == m*(p+1):
            if max(j*c,k*b) > max(n*c,p*b):
                return True
            elif max(j*c,k*b) == max(n*c,p*b):
                if j > n:
                    return True
                elif j == n:
                    if k > p:
                        return True
    return False

def left1(i,j,k):
    L = []
    if i > 0 or j > 0:
        L.append([i,j,k])
        L.append([-i,-j,k])
    if k > 0:
        L.append([i+1,-j-1,k-1])
        L.append([-i-1,j+1,k-1])
    return L

def left2(i,j,k):
    L = []
    if j > 0 or k > 0:
        L.append([i,j,k])
        L.append([i,-j,-k])
    if i > 0:
        L.append([i-1,j+1,-k-1])
        L.append([i-1,-j-1,k+1])
    return L

def right1(a,b,c,r,s,t):
    if (a-r > 0 and -b-s < 0) or (a-r < 0 and -b-s > 0):
        return [abs(a-r)-1,abs(-b-s)-1,t+1]
    return [abs(a-r),abs(-b-s),t]

def right2(a,b,c,r,s,t):
    if (-b-s > 0 and c-t < 0) or (-b-s < 0 and c-t > 0):
        return [r+1,abs(-b-s)-1,abs(c-t)-1]
    return [r,abs(-b-s),abs(c-t)]

def generatingset(a,b,c):
    GS = []
    MGS = []
    for i in range(a):
        for j in range(b):
            k = 0
            while b*k <= 2*c*(b-1):
                GS.append([i,j,k])
                k += 1
    for T in GS:
        rewrite = False
        for L in left1(T[0],T[1],T[2]):
            R = right1(a,b,c,L[0],L[1],L[2])
            if greaterthan(a,b,c,T[0],T[1],T[2],R[0],R[1],R[2]):
                rewrite = True or rewrite
        for L in left2(T[0],T[1],T[2]):
            R = right2(a,b,c,L[0],L[1],L[2])
            if greaterthan(a,b,c,T[0],T[1],T[2],R[0],R[1],R[2]):
                rewrite = True or rewrite
        if not rewrite:
            MGS.append(T)
    return MGS
\end{verbatim}

Using the code listed above, we obtain the generating sets listed in the introduction for $S(M(2,-2,2))$, $S(M(3,-2,3))$, and $S(M(3,-2,5))$, and we find that our generating set is minimal for $S(M(2,-2,2);\Z[A^{\pm 1}],A)$.  

As for observing minimality of our generating sets for $S(M(3,-2,3);R[A^{\pm 1}],A)$ and $S(M(3,-2,5);R[A^{\pm 1}],A)$, we might hope to consider $S(M(3,-2,3);R,-1)$ and $S(M(3,-2,5);R,-1)$, as they are isomorphic to the skein algebras of their fundamental groups, which are generated by representatives of conjugacy classes.  For $S(M(3,-2,3);R[A^{\pm 1}],A)$, however, this will not help, as only three of the conjugacy classes of the binary tetrahedral group are self-inversive:  $S(M(3,-2,3);R,-1)$ can be generated by five elements.  See Przytycki and Sikora \cite{PS}.

Still, for $S(M(3,-2,5);R[A^{\pm 1}],A)$, we can hope to gain some insight, as its conjugacy classes are self-inversive, and since we have the following result:

\begin{propo}
Suppose that a set $L = \{L_1, \ldots, L_n\}$  of links in $M$ represents a generating set for $S(M;R[A^{\pm 1}],A)$.
\begin{enumerate}

\item If $L$ yields a minimal generating set for $S(M;R,-1)$, then $L$ represents a minimal generating set for $S(M;R[A^{\pm 1}],A)$. 

\item If $L$ yields a linearly independent set for $S(M;R,-1)$ and $S(M;R[A^{\pm 1}],A)$ has no $(A+1)$ torsion, then $L$ represents a basis for $S(M;R[A^{\pm 1}],A)$. 

\item If $L$ yields a linearly independent set for $S(M;R,-1)$ and $S(M;R[A^{\pm 1}],A)$ has torsion, then $S(M;R[A^{\pm 1}],A)$ has $(A+1)$ torsion. 

\end{enumerate}
\label{prop:min}
\end{propo}

\begin{proof}
(1) Suppose that $L_n = f_1(A) L_1 + \cdots + f_{n-1}(A) L_{n-1}$ in $S(M;R[A^{\pm 1}],A)$.
Then in $S(M;R,-1)$, $L_n = f_1(-1) L_1 + \cdots + f_{n-1}(-1) L_{n-1}$, a contradiction.

(2) Suppose that $f_1(A) L_1 + \cdots + f_n(A) L_n = 0$ in $S(M;R[A^{\pm 1}],A)$.
Then in $S(M;R,-1)$, $f_1(-1) L_1 + \cdots + f_n(-1) L_n = 0$.
$L_1, \ldots, L_n$ is a basis of $S(M;R,-1)$, so $f_i(-1) = 0$ for each $i$, and thus $(A+1) | f_i$ for each $i$.
Hence, for some $g_1, \ldots g_n$, $(A+1)(g_1(A) L_1 + \cdots + g_n(A) L_n) = 0$.
$S(M;R[A^{\pm 1}],A)$ has no $(A+1)$ torsion, so $g_1(A) L_1 + \cdots + g_n(A) L_n = 0$.
Hence, $S(M;R[A^{\pm 1}],A)$ is free.

(3) If $L$ yields a linearly independent set for $S(M;R,-1)$ and $S(M;R[A^{\pm 1}],A)$ has torsion, then $L$ cannot represent a basis, and hence $S(M;R[A^{\pm 1}],A)$ must have $(A+1)$ torsion by (2).
\end{proof}

\begin{rem}
The existence of torsion is a topic of particular interest in skein theory.  For example, McLendon has studied $(A+1)$ torsion in \cite{M06}.
\end{rem}

Let $G$ be the binary icosahedral group, with presentation $\left< r,s | r^5 = s^3 = (rs)^2 \right>$.
Since $G$ is finite, the skein algebra of $G$ over $\C$ is isomorphic to $\C[X(G)]$, the $SL(2,\C)$ character variety of $G$ (\cite{PS}, see also Bullock \cite{Bu97_2}).
 
Let $\sigma_0$ be the trivial 2-dimensional representation of $G$, let $\sigma_1$ be the representation of $G$ that sends $r$ and $s$ to 

$$A_1 = \frac{1}{5}\left[
\begin{array}{cc}
-3e_5-e_5^2+e_5^3-2e_5^4 & e_5-3e_5^2-2e_5^3-e_5^4 \\
e_5+2e_5^2+3e_5^3-e_5^4 & -2e_5+e_5^2-e_5^3-3e_5^4
\end{array}\right]$$  
 
and

$$B_1 = \frac{1}{5}\left[
\begin{array}{cc}
-e_5-2e_5^2-3e_5^3-4e_5^4 & 2e_5-e_5^2+e_5^3-2e_5^4 \\ 
2e_5-e_5^2+e_5^3-2e_5^4 & -4e_5-3e_5^2-2e_5^3-e_5^4
\end{array}\right],$$

respectively, and let $\sigma_2$ be the representation of $G$ that sends $r$ and $s$ to   

$$A_2 = \left[
\begin{array}{cc}
e_5-e_5^2 & -e_5^2-e_5^4 \\
-e_5-e_5^4 & -e_5-e_5^3 
\end{array}\right] \mbox{ and } 
B_2 = \left[
\begin{array}{cc}
1 & -e_5^3 \\
e_5^2 & 0 
\end{array}\right],$$ respectively, where $e_5 = e^{\frac{2 \pi i}{5}}$.

Using GAP \cite{GAP4}, we can see that $\sigma_0$, $\sigma_1$, and $\sigma_2$ are $SL(2,\C)$ representations of $G$, and any $SL(2,\C)$ representation $\sigma$ of $G$ is equivalent to one of them:  if irreducible, $\sigma$ is equivalent to $\sigma_1$ or $\sigma_2$, and if reducible, $\sigma$ is equivalent to $\sigma_0$, since $G$ is perfect.  See Culler and Shalen \cite{CS}.

Let $\chi_0$, $\chi_1$, and $\chi_2$ be the characters of $\sigma_0$, $\sigma_1$, and $\sigma_2$, respectively, 
and for each $g \in G$, let $\tau_g$ be the evaluation map defined on the characters of $G$ by $\tau_g(\chi) = \chi(g)$.  Note that since $1, r, r^2, r^3, r^4, r^5, rs, s,$ and $s^2$ represent the conjugacy classes of $G$, $\C[X(G)]$ is generated by $\tau_1, \tau_r, \tau_{r^2}, \tau_{r^3}, \tau_{r^4}, \tau_{r^5}, \tau_{rs}, \tau_s,$ and $\tau_{s^2}$. 

$$\begin{array}{c|ccccccccc}
& \tau_1 & \tau_r & \tau_{r^2} & \tau_{r^3} & \tau_{r^4} & \tau_{r^5} & \tau_{rs} & \tau_s & \tau_{s^2} \\
\hline
\chi_0 & 2 & 2 & 2 & 2 & 2 & 2 & 2 & 2 & 2 \\ 
\chi_1 & 2 & -e_5-e_5^4 & e_5^2+e_5^3 & -e_5^2-e_5^3 & e_5+e_5^4 & -2 & 0 & 1 & -1 \\  
\chi_2 & 2 & -e_5^2-e_5^3 & e_5+e_5^4 & -e_5-e_5^4 & e_5^2+e_5^3 & -2 & 0 & 1 & -1  
\end{array}$$

From the table, we can see that the following relations hold in $\C[X(G)]$:

\begin{itemize}
\item $\tau_{s^2} = 3 \tau_s - 2 \tau_1$
\item $\tau_{rs} = 2 \tau_s - \tau_1$
\item $\tau_{r^5} = 4 \tau_s - 3 \tau_1$
\item $\tau_{r^4} = 4 \tau_s - \tau_r - 2 \tau_1$
\item $\tau_{r^3} = 3 \tau_s - \tau_r - \tau_1$
\item $\tau_{r^2} = \tau_s + \tau_r - \tau_1$
\end{itemize}

Furthermore, $\{\tau_1, \tau_r, \tau_s \}$ are linearly independent in $\C[X(G)]$, since the matrix

$$\left[
\begin{array}{ccc}
\tau_1(\chi_0) & \tau_r(\chi_0) & \tau_s(\chi_0) \\
\tau_1(\chi_1) & \tau_r(\chi_1) & \tau_s(\chi_1) \\
\tau_1(\chi_2) & \tau_r(\chi_2) & \tau_s(\chi_2) 
\end{array}\right]
=
\left[
\begin{array}{ccc}
2 & 2 & 2 \\ 
2 & -e_5-e_5^4 & 1 \\  
2 & -e_5^2-e_5^3 & 1   
\end{array}\right]$$

is invertible.  

Thus, $S(M(3,-2,5);\C,-1)$ is 3-dimensional, and so, we cannot use Proposition \ref{prop:min} to demonstrate that our generating set for $S(M(3,-2,5);\C[A^{\pm 1}],A)$ is minimal.  Hence, we are left with the following:

\begin{ques}  Is $\{1, z, z^2, z^3, z^4, z^5, y, x, x^2\}$ a minimal generating set for $\\S(M(3,-2,5);R,A)$ for some ring $R$ and unit $A$?  If not, is $S(M(3,-2,5);R,A)$ generated by $\{1, z, x\}$ for every ring $R$ and unit $A$?
\end{ques}

\section{Acknowledgements}
\label{sec:Acknowledgements}

The author would like to thank Patrick Gilmer for his suggestions and encouragement.

\end{document}